\documentclass[a4paper,11pt]{article}

\usepackage{latexsym}
\usepackage{amssymb}
\usepackage{amsthm}
\usepackage{amsmath}
\usepackage{amssymb}
\usepackage{enumerate}
\usepackage{color}
\usepackage{graphicx}

\def\endproof{\hfill$\Box$}

\newcommand{\eps}{\varepsilon}

\newcommand{\dd}{\hspace{1pt}{\rm d}\hspace{0.0pt}}
\newcommand{\ee}{{\rm e}\hspace{1pt}}

\newcommand{\R}{\mathbb R}
\newcommand{\N}{\mathbb N}

\newcommand{\cB}{{\mathcal B}}

\newcommand{\cD}{{\mathcal D}}
\newcommand{\cE}{{\mathcal E}}
\newcommand{\cF}{{\mathcal F}}

\newcommand{\cH}{{\mathcal H}}

\newcommand{\cK}{{\mathcal K}}
\newcommand{\cL}{{\mathcal L}}
\newcommand{\cM}{{\mathcal M}}
\newcommand{\cN}{{\mathcal N}}

\newcommand{\cP}{{\mathcal P}}
\newcommand{\cR}{{\mathcal R}}
\newcommand{\cS}{{\mathcal S}}

\newcommand{\rE}{\ensuremath{\mathrm{E}}}
\newcommand{\rV}{\ensuremath{\mathrm{V}}}

\newcommand{\inv}{^{-1}}

\newtheorem{thm}{Theorem}[section]
\newtheorem{lemma}[thm]{Lemma}
\newtheorem{prop}[thm]{Proposition}
\newtheorem{defn}[thm]{Definition}
\newtheorem{cor}[thm]{Corollary}
\newtheorem{rem}[thm]{Remark}
\newtheorem{ex}[thm]{Example}

\title{Generalized stochastic processes revisited}

\vspace{5mm}
\author{Michael Oberguggenberger\footnote{Unit for Engineering Mathematics, University of Innsbruck, Technikerstr.\ 13, A-6020 Innsbruck, Austria, michael.oberguggenberger@uibk.ac.at}
         }
\date{}

\begin{document}

\maketitle

{\bf Abstract.}
The paper addresses the question whether a random functional, a map from a set $E$ into the space of real-valued measurable functions on a probability space, has a measurable version with values in $\R^E$.
Similarly, one may ask whether linear random functionals have versions in the algebraic dual. Most importantly, it can be asked which locally convex topological vector spaces $E$ have the ``regularity property'' that any linear random functional on $E$ has a version with values in the dual $E'$, an important issue in the theory of generalized stochastic processes. It has been shown by It\^{o} and Nawata that this is the case when $E$ is nuclear. However, the question of uniqueness has only been partially answered. We build up a framework where these and related questions can be clarified in terms of spaces and mappings. We study classes of spaces $E$ (beyond nuclear spaces) with the said regularity property, prove a seemingly new uniqueness result and exhibit various examples and counterexamples.
\vspace{5mm}
\vspace{5mm}
%
%
\section{Introduction}
\label{Sec:intro}
%
%
This paper takes up the question of relating different notions of generalized stochastic processes to each other, mainly motivated by processes with values in the space of Schwartz distributions $\cD'(\R^n)$.
Gel'fand and Vilenkin \cite{GelfandVilenkin:1964} (1961/1964) considered such ``generalized random functions'' as families of random variables indexed by test functions $\varphi\in \cD(\R^n)$, which are linear in $\varphi$ and continuous in distribution
on $\cD(\R^n)$. On the other hand, Ullrich \cite{Ullrich:1957} (1956) introduced ``random Schwartz distributions'' as weakly measurable maps on a probability space with values in $\cD'(\R^n)$; this notion has been widely adopted (e.g. \cite{Hida:1980}).
Further, ``random distributions'' -- continuous linear functionals on $\cD(\R^n)$ with values in a Hilbert space of second order random variables -- were introduced by It\^{o} \cite{Ito:1954} (1953).

It has been established by It\^{o} and Nawata in \cite{Ito:1984,ItoNawata:1983} that any generalized random function in the sense of Gel'fand-Vilenkin has a version which is a random Schwartz distribution in the sense of Ullrich. This nontrivial result rests on the nuclearity of $\cD(\R^n)$.
Also, the linearity of generalized random functionals implies that continuity in distribution is equivalent with continuity in probability. This shows that It\^{o} linear functionals also belong to the Gel'fand-Vilenkin class.

The aim of the paper is to systematically study the relations between these classes of stochastic processes in terms of spaces and mappings, with and without linearity and/or continuity.

Let $(\Omega,\Sigma,P)$ be a probability space and $E$ a set. A \emph{random functional} \cite[Section 2.3]{Ito:1984} is a map $Z:E\times\Omega \to \R$ such that the map $\omega\to Z(e,\omega)$ is measurable for every $e\in E$.
The random functional $Z$ can be viewed from two sides: either as a map $X:\Omega\to\R^E$ or as a map $Y:E\to\R^\Omega$. This dichotomy underlies the difference in interpretation of Ullrich vs. Gel'fand-Vilenkin.
We start by observing that the map
\[
   \kappa:L^0(\Omega\!:\!\R^E) \to (L^0(\Omega))^E, \qquad \kappa X(e)(\omega) = X(\omega)(e)  
\]
is well-defined, where $L^0(\Omega\!:\! F)$ denotes the space of (equivalence classes of) measurable maps from $(\Omega,\Sigma)$ into a measurable space $(F,\cF)$. It turns out that the map $\kappa$ is always surjective, but not necessarily injective. 
On the next level, we let $E$ be a linear space with algebraic dual $E^\ast$ and consider the restriction of $\kappa$ as
\[
   \kappa_0:L^0(\Omega\!:\!E^\ast) \to \cL(E\!:\!L^0(\Omega))
\]
and study its surjectivity/injectivity, where $\cL(E\!:\!L^0(\Omega))$ denotes the space of linear maps from $E$ to $L^0(\Omega)$. The central issue is adding continuity, yielding the map
\[
   \overline{\kappa}: L^0(\Omega\!:\! E') \to \cL_{\rm c}(E\!:\!L^0(\Omega))
\]
where $E'$ denotes the sequentially continuous dual of $E$ and $\cL_{\rm c}(E\!:\!L^0(\Omega))$ the space of sequentially continuous linear maps from $E$ to $L^0(\Omega)$.
A more complete diagram relating the mentioned spaces can be found in Section\;\ref{Sec:contlinrandfun}.

Let us temporarily say that a locally convex topological vector space has the \emph{surjectivity property}, if the map $\overline{\kappa}$ is surjective for whatever probability space $(\Omega,\Sigma, P)$.
The mentioned result of It\^{o} and Nawata says that every nuclear space $E$ has the surjectivity property. We exhibit several examples where the surjectivity property fails when $E$ is not nuclear.
We also investigate whether the surjectivity property is hereditary. For example, we show that if $E = \sum_{\gamma\in\Gamma}E_\gamma$ is the locally convex direct sum of spaces $E_\gamma$ with the surjectivity property, then $E$ has the surjectivity property (and conversely). If $\Gamma$ is uncountable and each $E_\gamma$ is nuclear, $E$ provides an example of a space with the surjectivity property which is not nuclear.

The result of It\^{o} and Nawata also gives a certain uniqueness of the regular version of a process $\cL_{\rm c}(E\!:\!L^0(\Omega))$ in $L^0(\Omega\!:\! E')$, but falls short of proving injectivity of the map $\overline{\kappa}$.
We show that $\overline{\kappa}$ is injective when $E$ is sequentially separable. We also exhibit various example of the failure of injectivity. A crucial observation is that equivalence of measurable functions 
in $L^0(\Omega\!:\! E')$ is stronger than equivalence of function values in $\cL_{\rm c}(E\!:\!L^0(\Omega))$.

The results of the paper provide some progress on the ``enticing open question'' posed by Pilipovi\'{c} and Sele\v{s}i in \cite[p. 252]{PilipovicSelesi:2007}.

The plan of the paper is as follows. In Section\;\ref{Sec:motivation} we describe the initial question of the paper in more detail and introduce some notation. Random functionals are addressed in Section\;\ref{Sec:randfun},
linear random functionals in Section\;\ref{Sec:linrandfun}. The central parts are Sections\;\ref{Sec:contlinrandfun} and \ref{Sec:regularity}. In Section\;\ref{Sec:contlinrandfun}, continuous linear random functionals are considered. Various notions of continuity are related to each other, and the relations of the resulting spaces are worked out. The role of sequential separability is also pinned down. 

Section\;\ref{Sec:regularity} starts out by recalling the central result of It\^{o} and Nawata that any continuous linear random functional on a nuclear space $E$ has a version in $E'$. Combining it with our injectivity result, we obtain that 
$\overline{\kappa}$ is a bijection when $E$ is a sequentially separable nuclear space.
It is also shown, in particular, that the surjectivity property is inherited by locally convex direct sums, complemented subspaces and by countable regular inductive limits.
Section\;\ref{Sec:counter} provides various counterexamples, in particular, simple examples showing that the theorem of It\^{o} and Nawata does not hold without the nuclearity assumption, in general. The mapping properties between the various spaces considered in Section\;\ref{Sec:contlinrandfun} are further analyzed. Finally, the Appendix serves to record the theorem of It\^{o} and Nawata, which originally slightly goes beyond nuclearity of the space $E$ itself. We also record the version of Walsh \cite{Walsh:1984} who proves the surjectivity property for complete, countably Hilbertian nuclear spaces.

We end the introduction with a few references to the literature. The issue whether a linear random functional (of Gel'fand-Vilenkin-type) on some topological vector space has a version (of Ullrich-type) in the continuous dual has been pointed out
in a number of papers, e.g. \cite{Brzezniak:2015}, \cite{Mirkov:2009}, who refer to the simplified special version of It\^{o}'s proof in Walsh \cite[Corollary 4.2]{Walsh:1984}.
A detailed comparison of It\^{o}-type random distributions with Ullrich-type functionals can be found in \cite{PilipovicSelesi:2007}, together with many examples and counterexamples,
mainly on $\cK\{M_p\}$-spaces (see also \cite{Swartz:1971} for this setting).

Generalized stochastic processes can also be defined as measures on the dual space. This topic, together with characteristic functionals, remains outside the scope of this paper.
The relation is extensively discussed in It\^{o}'s monograph \cite{Ito:1984} and at many other places, e.g. \cite{Hida:1970}. Here the question whether such a measure on the algebraic dual gives probability one to the continuous dual
(and then allows one to construct Ullrich-type-versions of the corresponding process) has been extensively discussed in \cite{Dudley:1969}.
Similar questions arise with the earlier concept of ``weak distributions'' in the sense of Segal \cite{Segal:1954}, introduced as random variables on the continuous dual; see also \cite{Getoor:1957}, \cite{Gross:1963}.

\section{Motivation and notations}
\label{Sec:motivation}

In order to motivate the exposition, we first formulate the concepts of Gel'fand-Vilenkin, Ullrich and It\^{o} in more detail; for the notation of spaces of test functions and Schwartz distributions, see e.g. \cite{Horvath:1966}.

Gel'fand and Vilenkin define a ``generalized random function'' in \cite[Section III.1.2.]{GelfandVilenkin:1964} as a stochastic process $X$ indexed by $\cD = \cD(\R^n)$ such that

(i) For all $\varphi\in\cD$, $X(\varphi)$ is an $\R$-valued random variable;

(ii) $X$ is linear in the sense that for all $\varphi,\psi\in\cD$, $\alpha,\beta\in\R$, $X(\alpha\varphi + \beta\psi)$ has the same distribution as $\alpha X(\varphi) + \beta X(\psi)$;

(iii) if $\varphi_n \to \varphi$ in $\cD$, then $X(\varphi_n) \to X(\varphi)$ in distribution.

Ullrich \cite{Ullrich:1957} starts by introducing the $\sigma$-algebra on $\cD'$ generated by the sets $\{T\in\cD': T(\varphi) < c\}$, $\varphi\in\cD$, $c\in\R$. Given a measurable space $(\Omega,\Sigma)$,
a ``random Schwartz distribution'' is a measurable transformation $X:\Omega\to \cD'$. If $(\Omega,\Sigma,P)$ is a probability space, this is the same as saying that
$X$ is a family of random variables $\{X(\varphi,\omega):\varphi\in\cD\}$ such that $X(\cdot,\omega)\in\cD'$ for almost all $\omega\in\Omega$, see Hida \cite{Hida:1980}, last paragraph in Section 1.3.(ii).

Finally, It\^{o} \cite{Ito:1954} introduces ``random distributions'' as continuous linear functionals on $\cD$ with values in the Hilbert space $\cH$ of random variables
with mean zero and finite variance.

We set out to analyze the relations between these definitions, replacing $\cD$ by $E$, where we will successively take for $E$ a set, a linear space, and a locally convex topological vector space.

In the sequel, $(\Omega,\Sigma,P)$ will denote a probability space. Given a measurable space $(F,\cF)$, the space $L^0(\Omega\!:\!F)$ denotes the set of (equivalence classes) of measurable maps $(\Omega,\Sigma)\to(F,\cF)$.
In case $F=\R$, we will take $\cF = \cB(\R)$, the Borel $\sigma$-agebra on $\R$ and simply write $L^0(\Omega)$ for $L^0(\Omega\!:\!\R)$.
We will always assume that all representatives of measurable maps are defined on all of $\Omega$ and are finite-valued.

Let $E$ be a set. The space $\R^E$ of maps $E\to\R$ will be equipped with the Kolmogorov $\sigma$-algebra
$\cB_K(\R^E)$, which is generated by the sets
\[
   \pi\inv_e(B), B \in \cB(\R), e\in E,
\]
where $\pi_e$ denotes the projection on the e-th slot. (This $\sigma$-algebra is often referred to as the product $\sigma$-algebra, but it is different from the $\sigma$-algebra generated by the product topology, unless $E$ is countable.)

When $G,H$ are sets, we shall occasionally use the notation $\cM(G\!:\!H)$ for the family of all maps $G\to H$, that is, $\cM(G\!:\!H)= H^G$.

\section{Random functionals}
\label{Sec:randfun}

A \emph{random functional} on a set $E$ is a map $Z:E\times\Omega \to \R$ such that the map $\omega\to Z(e,\omega)$ is measurable for every $e\in E$.
As noted in the Introduction, the random functional $Z$ can be viewed from two sides: either as a map $X:\Omega\to\R^E$ or as a map $Y:E\to\R^\Omega$. 
We will build up a framework yielding a canonical linear map that allows one to connect the two viewpoints via
\[
   \kappa:L^0(\Omega\!:\!\R^E) \to (L^0(\Omega))^E. 
\]
We begin by defining the map $\kappa$ on single representatives as
\[
   \kappa:\cM(\Omega\!:\!\R^E) \to (\cM(\Omega\!:\!\R))^E, \quad \kappa X(e)(\omega) = X(\omega)(e).
\]
\begin{lemma}\label{lem:measurability}
$X:\Omega\to \R^E$ is $(\Omega,\Sigma)-(\R^E,\cB_K(\R^E))$-measurable if and only if $\kappa X(e)$ is $(\Omega,\Sigma)-(\R,\cB(\R))$-measurable for all $e\in E$.
\end{lemma}
\emph{Proof.}
Since the sets $\pi\inv_e(B)$, $B\in \cB(\R)$, $e\in E$ generate the Kolmogorov $\sigma$-alebra, $X$ is measurable if and only if
all sets $\{\omega\in\Omega:X(\omega)\in \pi\inv_e(B)\}$ ($B\in\cB(\R)$, $e\in E$) belong to $\Sigma$. On the other hand, $\kappa X(e)$ is measurable if and only if
all sets $\{\omega\in\Omega: \kappa X(e)(\omega) \in B\}$ ($B\in\cB(\R)$) belong to $\Sigma$. But
\begin{eqnarray*}
 & &\{\omega\in\Omega: \kappa X(e)(\omega) \in B\} = \{\omega\in\Omega: X(\omega)(e) \in B\}\\
 & & = \{\omega\in\Omega: \pi_eX(\omega) \in B\} = \{\omega\in\Omega: X(\omega)\in \pi\inv_e(B)\},
\end{eqnarray*}
so the assertion follows. \endproof

In the next step, we go over to equivalence classes. Let $X_1,X_2$ be measurable maps $\Omega \to \R^E$. Then $X_1$ and $X_2$ are strictly equivalent
\cite[p. 258]{ItoNawata:1983} (or strict versions of each other), if and only if
\[
\exists \Omega_1\subset \Omega, P(\Omega_1) = 1 \quad {\rm such\ that}\quad \forall \omega\in \Omega_1, X_1(\omega) = X_2(\omega).
\]
The latter assertion means that for all $\omega\in \Omega_1$ and all $e\in E$, $X_1(\omega)(e) = X_2(\omega)(e)$. On the other hand, if $Y_1,Y_2$ are maps from $E$ into $\cM(\Omega\!:\!\R)$ such that
$Y_1(e)$ and $Y_2(e)$ are measurable for each $e\in E$, equivalence takes place in $L^0(\Omega)$, meaning that
\[
\forall e\in E \ \exists \Omega_{1,e}\subset \Omega, P(\Omega_{1,e}) = 1 \quad {\rm such\ that}\quad \forall \omega\in \Omega_{1,e}, Y_1(e)(\omega) = Y_2(e)(\omega).
\]
In the language of stochastic processes indexed by $E$, this means that $Y_1$ is a version of $Y_2$.

\begin{lemma}\label{lem:strict}
If $X_1$ is a strict version of $X_2$, then $\kappa X_1$ is a version of $\kappa X_2$.
\end{lemma}
\emph{Proof.}
By assumption, there is a subset $\Omega_1$ of full measure such that for all $\omega\in \Omega_1$ and all $e\in E$, $X_1(\omega)(e) = X_2(\omega)(e)$. Since $\kappa X_i(e)(\omega) = X_i(\omega)(e)$, this in turn implies the equivalence of $\kappa X_1(e)$ and $\kappa X_2(e)$ as measurable maps $\Omega\to \R$, with $\Omega_{1,e} = \Omega_1$ for all $e$. \endproof
\begin{cor}\label{cor:kappasurjective}
The map $\kappa$ does not depend on the choice of representatives, thus can be extended to a map
\[
 \kappa:L^0(\Omega\!:\!\R^E) \to (L^0(\Omega))^E,
\]
and $\kappa$ is surjective.
\end{cor}
\emph{Proof.}
Denoting equivalence classes of measurable maps by square brackets, the assignment
\[
  [X] \to \kappa[X] = [\kappa X(e)]_{e\in E}
\]
is well-defined by Lemma\;\ref{lem:strict}. To prove surjectivity, let $Y\in (L^0(\Omega))^E$. For each $e\in E$, $Y(e)$ is an equivalence class of measurable functions. Invoking the axiom of choice, we may choose a representative $\overline{Y(e)}\in Y(e)$ for each $e\in E$. Define $X\in \cM(\Omega\!:\!\R^E)$ by
\[
   X(\omega)(e) = \overline{Y(e)}(\omega).
\]
Then $\kappa X(e)(\omega) = X(\omega)(e) = \overline{Y(e)}(\omega)$ is $(\Omega,\Sigma)-(\R,\cB(\R))$-measurable. By Lemma\;\ref{lem:measurability}, $X$ is
$(\Omega,\Sigma)-(\R^E,\cB_K(\R^E))$-measurable, and $\kappa X(e) = \overline{Y(e)} \in Y(e)$. Thus $[X]\in L^0(\Omega\!:\! \R^E)$, and
$\kappa[X] = Y$.    \endproof

The map $\kappa$ is not necessarily injective (Example\;\ref{rem:kappanot1-1}).

From now on, unless clarity requires it, we will suppress the square bracket notation, and -- as is common practice -- employ the same notation for a representative and its class.

\section{Linear random functionals}
\label{Sec:linrandfun}

In this section, $E$ will be a vector space over $\R$, and $E^\ast$ its algebraic dual. The algebraic dual is equipped with the restriction of the Kolmogorov $\sigma$-algebra $E^\ast \cap \cB_K(\R^E)$. Thus a map $\Omega\to E^\ast$ is measurable, iff $\pi_eX$ belongs to $L^0(\Omega)$ for all $e\in E$. The space of linear mappings $E\to L^0(\Omega)$ is denoted by $\cL(E\!:\!L^0(\Omega))$.

Now $X\in L^0(\Omega\!:\!E^\ast)$ means that there is a set $\Omega_1\subset \Omega$ of full measure, such that $X(\omega)$ is linear in $E$. That is,
\begin{eqnarray*}
   \exists \Omega_1, P(\Omega_1 = 1): \forall \omega\in\Omega_1, \alpha,\beta\in\R, e,f\in E,\\
    X(\omega)(\alpha f + \beta g) = \alpha X(\omega)(e) + \beta X(\omega)(f).
\end{eqnarray*}
On the other hand, $Y\in \cL(E\!:\!L^0(\Omega))$ means that Y is linear in $E$ with values in $L^0(\Omega)$. That is,
\begin{eqnarray*}
\forall \alpha,\beta\in\R, e,f\in E, \exists \Omega_1,  P(\Omega_1 = 1): \forall \omega\in\Omega_1,\\
    X(\omega)(\alpha f + \beta g) = \alpha X(\omega)(e) + \beta X(\omega)(f).
\end{eqnarray*}
To be sure, in the second case $\Omega_1$ may depend on $\alpha,\beta,e,f$, while in the first case it does not.

If $X\in L^0(\Omega\!:\!E^\ast)$, then clearly $\kappa X\in \cL(E\!:\!L^0(\Omega))$. Thus the restriction of $\kappa$ to $L^0(\Omega\!:\!E^\ast)$ results in a map
\[
   \kappa_0: L^0(\Omega\!:\!E^\ast) \to \cL(E\!:\!L^0(\Omega)),
\]
which is also well-defined.

\begin{prop}\label{prop:kappa0surjective}
The map $\kappa_0: L^0(\Omega\!:\!E^\ast) \to \cL(E\!:\!L^0(\Omega))$ is surjective.
\end{prop}
\emph{Proof.}
Let $Y\in\cL(E\!:\!L^0(\Omega))$. Choose an  an algebraic basis $(e_\lambda)_{\lambda\in\Lambda}$ of $E$. 
Let
\[
   \R^{(\Lambda)} = \{f \in \R^\Lambda: f(\lambda)\neq 0\ {\rm for\ at\ most\ finitely\ many\ } \lambda\in\Lambda\}.
\]
Then $E$ is isomorphic with $\R^{(\lambda)}$; the isomorphism $\theta$ sends $e_\lambda$ to the $\lambda$th unit vector in $\R^{(\Lambda)}$. Its transpose $\theta'$ is an isomorphism between $\R^\Lambda$ and $E^\ast$.
Choose representatives $\overline{Y(e_\lambda)}\in Y(e_\lambda)$ and set
\[
   X(\omega) = \theta'\Big(\big(\overline{Y(e_\lambda)}(\omega)\big)_{\lambda\in\Lambda}\Big).
\]
Clearly $X(\omega) \in E^\ast$ for all $\omega\in\Omega$. If $e = \sum_{\rm finite}\alpha_\lambda e_\lambda\in E$, then $X(\omega)(e) = \sum_{\rm finite}\alpha_\lambda \overline{Y(e_\lambda)}(\omega)$, so
$X:\Omega\to E^\ast$ is measurable. Finally, 
\[
\big[\sum_{\rm finite}\alpha_\lambda \overline{Y(e_\lambda)}\big] = \sum_{\rm finite}\alpha_\lambda \big[\overline{Y(e_\lambda)}] = \sum_{\rm finite}\alpha_\lambda Y(e_\lambda) = Y(e),
\]
since $Y$ is linear. This shows that $[\kappa_0 X(e)] = Y(e)$ for all $e\in E$ and so $\kappa_0[X] = Y$.
\endproof

The surjectivity of $\kappa_0$ can also be obtained as a corollary to Proposition\;\ref{prop:directsum} (Remark\;\ref{rem:finestlctopology}).
The map $\kappa_0$ need not be injective (Example\;\ref{rem:kappa0not1-1}). However, we have the following positive result.
\begin{prop}\label{prop:countable}
Suppose the dimension of $E$ is at most countable. Then $\kappa_0: L^0(\Omega\!:\!E^\ast) \to \cL(E\!:\!L^0(\Omega))$ is injective.
\end{prop}
\emph{Proof.}
An element $X \in L^0(\Omega\!:\!E^\ast)$ is determined by its action on an algebraic basis $\{e_\iota:\iota\in I\}$.
Indeed, if $e\in E$ and $e = \sum_{\rm finite} \lambda_\iota e_\iota$, then
\[
   X(\omega)(e) = \sum_{\rm finite} \lambda_\iota X(\omega)(e_\iota), \qquad \kappa_0 X(e)(\omega) = \sum_{\rm finite} \lambda_\iota \kappa_0 X(e_\iota)(\omega).
\]
Now $\kappa_0 X(e) = 0$ in $L^0(\Omega)$ for all $e\in E$ if and only if, for all $\iota\in I$, $\kappa_0 X(e_\iota)(\omega) = 0$ on sets $\Omega_\iota$ of full measure. But then $X(\omega)(e_\iota) = 0$ for all $\iota\in I$
on $\bigcap_{\iota\in I} \Omega_\iota$, which also has full measure. Thus $X = 0$ in $L^0(\Omega\!:\!E^\ast)$. \endproof

\section{Continuous linear random functionals}
\label{Sec:contlinrandfun}

In what follows, the term \emph{locally convex space} will always mean a separated (i.e., Hausdorff) locally convex topological vector space. Let $E$ be a locally convex space. The space of \emph{sequentially continuous linear maps} $E\to\R$ will be denoted by $E'$. It coincides with the continuous dual of $E$, when $E$ is a \emph{Mazur space} \cite[Definition 8-6-3]{Wilansky:1978} (this happens in all application cases we have in mind). If $F$ is a topological vector space, we denote by $\cL_c(E\!:\!F)$ the sequentially continuous linear maps from $E$ to $F$.

Let $X \in \cL(E\!:\!L^0(\Omega))$ and let $e_n\to e$ in $E$. We can consider various notions of continuity of $X$:
\begin{list}{\leftmargin}{}
\item[(dist)] $X(e_n) \to X(e)$ in distribution. This is the formulation of Gel'fand-Vilenkin \cite{GelfandVilenkin:1964}.

\item[(prob)] $X(e_n) \to X(e)$ in probability. This is assumed e.g. in It\^{o} \cite{Ito:1984} and Walsh \cite{Walsh:1984}.

\item[(path)] $X(e_n) \to X(e)$ almost surely (= pathwise almost everywhere), i.e., $X(e_n)(\omega) \to X(e)(\omega)$ for all $\omega$ in a set of full measure (which may depend on $e_n$ and $e$).
\end{list}

The subspaces of $\cL(E\!:\!L^0(\Omega))$, equipped with the respective convergence, will be denoted by $\cL(E\!:\!L^0(\Omega))_{\rm dist}$,
$\cL(E\!:\!L^0(\Omega))_{\rm prob}$, and $\cL(E\!:\!L^0(\Omega))_{\rm path}$.

\begin{rem}\label{rem:convergences}
(a) If $X\in\cL(E\!:\!L^0(\Omega))$ then $X(e_n)\to X(e)$ iff $X(e_n - e)\to 0$, by linearity. But convergence in distribution to a fixed constant is equivalent with
convergence in probability \cite[Section 3, p. 27]{Billingsley:1999}. Thus $\cL(E\!:\!L^0(\Omega))_{\rm dist} = \cL(E\!:\!L^0(\Omega))_{\rm prob}$.

(b) $L^0(\Omega)$ is a complete metrizable topological vector space with distance $d(f,g) = \rE(|f-g| \wedge 1)$. Convergence in metric is equivalent with convergence in probability
\cite[Section 9.3, Remark]{Loeve:1977}. Thus $\cL(E\!:\!L^0(\Omega))_{\rm prob} = \cL_{\rm c}(E\!:\!L^0(\Omega))$, the continuous linear maps from $E$ to $L^0(\Omega)$.

(c) For $p>0$, the imbedding $L^p(\Omega)\to L^0(\Omega)$ is continuous. This is obvious for $p\geq 1$; for $0 < p <1$, see \cite[Section 9.3, Application (2)]{Loeve:1977}.   

(d) Almost sure convergence implies convergence in probability.
\end{rem}
By (a) and (d) of Remark\;\ref{rem:convergences}, we have the following chain of inclusions:
\[
\cL(E\!:\!L^0(\Omega))_{\rm path} \subset  \cL(E\!:\!L^0(\Omega))_{\rm prob} = \cL(E\!:\!L^0(\Omega))_{\rm dist} \subset \cL(E\!:\!L^0(\Omega)).
\]
Further, by (b),
\[
\cL(E\!:\!L^0(\Omega))_{\rm prob} =  \cL_{\rm c}(E\!:\!L^0(\Omega))
\]
and, by (c),
\[
\cL_{\rm c}(E\!:\!L^p(\Omega)) \subset \cL_{\rm c}(E\!:\!L^0(\Omega))\quad{\rm for\ } 0\leq p\leq\infty.
\]
Similarly, for measurable maps with values in $E^\ast$, we have
\[
L^0(\Omega\!:\!E^\ast)_{\rm path} \subset  L^0(\Omega\!:\!E^\ast)_{\rm prob} = L^0(\Omega\!:\!E^\ast)_{\rm dist} \subset L^0(\Omega\!:\!E^\ast).
\]
Stipulating that the exceptional set arising in the definition of pathwise convergence in $L^0(\Omega\!:\!E^\ast)_{\rm path}$ does not depend on $e_n$ and $e$, we have
\[
   L^0(\Omega\!:\!E^\ast)_{\rm path} = L^0(\Omega\!:\!E').
\]
The relations can be summarized in the following diagram:
\[
\begin{array}{lllllll}
&&L^0(\Omega\!:\!\R^E) & \stackrel{\kappa}{\longrightarrow} &  \cM(E\!:\!L^0(\Omega)) &&\\[4pt]
&&\qquad\rotatebox[origin = c]{90}{$\hookrightarrow$} && \qquad\rotatebox[origin = c]{90}{$\hookrightarrow$} &&\\[4pt]
&&L^0(\Omega\!:\!E^\ast) &\stackrel{\kappa_0}{\longrightarrow} &  \cL(E\!:\!L^0(\Omega)) &&\\[4pt]
&&\quad\ \ \rotatebox[origin = c]{90}{$\hookrightarrow$} \ \iota &&\quad\ \ \rotatebox[origin = c]{90}{$\hookrightarrow$} \ \iota_0  &&\\[4pt]
&&L^0(\Omega\!:\!E^\ast)_{\rm dist} &\stackrel{\kappa_1}{\longrightarrow}&  \cL(E\!:\!L^0(\Omega))_{\rm dist}&&\\[4pt]
&&\qquad\rotatebox[origin = c]{90}{$=$} && \qquad\rotatebox[origin = c]{90}{$=$} &&\\[4pt]
&&L^0(\Omega\!:\!E^\ast)_{\rm prob} &\stackrel{\kappa_1}{\longrightarrow}& \cL(E\!:\!L^0(\Omega))_{\rm prob}  & = & \cL_{\rm c}(E\!:\!L^0(\Omega))\\[4pt]
&&\quad\ \ \rotatebox[origin = c]{90}{$\hookrightarrow$} \ \iota_1 &\!  {}^{\bar{\kappa}}\!{\nearrow} & \quad\ \ \rotatebox[origin = c]{90}{$\hookrightarrow$}\ \iota_2 &&\\[4pt]
L^0(\Omega\!:\!E')& =& L^0(\Omega\!:\!E^\ast)_{\rm path}  &\stackrel{\kappa_2}{\longrightarrow} & \cL(E\!:\!L^0(\Omega))_{\rm path} &&
\end{array}
\]
Here $\kappa_i$ denotes the restriction of $\kappa$ to the respective level, while $\bar{\kappa} = \kappa_1\circ\iota_1 = \iota_2\circ\kappa_2$.
All vertical arrows are injections.
The entry $\cL(E\!:\!L^0(\Omega))_{\rm dist}$ corresponds to the Gel'fand-Vilenkin
definition, $L^0(\Omega\!:\!E')$ corresponds to Ullrich and Hida. The spaces $\cL_{\rm c}(E\!:\!L^p(\Omega))$ (continuous linear maps on $E$ with values in $L^p(\Omega)$ enter on the right-hand side of the next to last line as subspaces of $\cL_{\rm c}(E\!:\!L^0(\Omega))$.

The central issue of this paper is the question of ``regularity'': Under what conditions does a linear random function belonging to $\cL_{\rm c}(E\!:\!L^0(\Omega))$ have a version with values in $E'$ -- or even an almost surely unique version?

\begin{defn}
\label{def:regular}
A locally convex space $E$ is said to have property ($\cR$), if for whatever probability space $(\Omega,\Sigma,P)$, the map
$\bar{\kappa}: L^0(\Omega\!:\!E')\to\cL_{\rm c}(E\!:\!L^0(\Omega))$ is bijective.
\end{defn}

Surjectivity means that for any $Y\in \cL_{\rm c}(E\!:\!L^0(\Omega))$ there is $X\in L^0(\Omega\!:\!E')$ such that $\bar{\kappa}X = Y$. Thus there is $\Omega_1\subset\Omega$ of full measure such that $X(\omega)\in E'$ for all $\omega\in \Omega_1$ and, for all $e\in E$, there is $\Omega_{1,e} \subset\Omega$ of full measure such that $X(\omega)(e) = Y(e)(\omega)$ for all $\omega\in \Omega_{1,e}$. In words, $X$ is a version of $Y$ with values in $E'$.

Note that one can assume without loss of generality that $\Omega_1 = \Omega$, setting e.g. $X(\omega) = 0$ for $\omega\not\in \Omega_1$.

Injectivity means that if $X,X'\in L^0(\Omega\!:\!E')$ and $\bar{\kappa}X = \bar{\kappa}X'$, then $X$ and $X'$ are strict versions of each other. In words: If for all $e\in E$ there is $\Omega_{1,e} \subset\Omega$ of full measure such that $X(\omega)(e) = X'(\omega)(e)$ for all $\omega\in \Omega_{1,e}$, then there is $\Omega_1\subset\Omega$ of full measure such that $X(\omega) = X'(\omega)$ in $E'$.

Positive examples of spaces with property ($\cR$) and permanence properties will be given in Section\;\ref{Sec:regularity}, counterexamples are in Section\;\ref{Sec:counter}.

\begin{rem}\label{rem:consequencesOfR}
(a) If the space $E$ has property ($\cR)$ then one has the relations
\[
  \cL_{\rm c}(E\!:\!L^p(\Omega)) \subset \cL_{\rm c}(E\!:\!L^0(\Omega)) = L^0(\Omega\!:\!E').
\]
(b) If the space $E$ has property ($\cR)$ then the map $\kappa_1$ is surjective, and the maps $\kappa_2$ and $\iota_2$ are both bijective.
Indeed, the surjectivity of $\bar{\kappa} = \kappa_1\circ\iota_1 = \iota_2\circ\kappa_2$ implies the surjectivity of $\kappa_1$ and of $\iota_2$; the injectivity of $\bar{\kappa}$ implies the injectivity of $\kappa_2$. But $\iota_2$ is injective by definition, thus it becomes a bijection, and then the same holds for $\kappa_2$.

It is not clear whether $\kappa_1$ is injective in this case (which in this situation is equivalent with the surjectivity of $\iota_1$).
\end{rem}

As will be observed in Section\;\ref{Sec:counter}, the maps $\kappa_0$, $\kappa_1, \kappa_2$ and $\bar{\kappa}$ need not be injective, in general. However, there is a convenient positive result.

\begin{prop}\label{prop:injective}
If $E$ is a sequentially separable locally convex space, then $\bar{\kappa}$ and $\kappa_2$ are injective.
\end{prop}
\emph{Proof.}
Let $X,X'\in L^0(\Omega\!:\!E')$ such that $\bar{\kappa} X = \bar{\kappa} X'$. For every $e\in E$, there is $\Omega_{1,e}\subset \Omega$ of full measure, such that
$X(\omega)(e) = X'(\omega)(e)$ for all $\omega\in\Omega_{1,e}$. Let $(e_n)_{n\in\N}$ be a countable dense subset of $E$. Set
\[
  \Omega_1 = \bigcap_{n=1}^\infty \Omega_{1,e_n}.
\]
Then $X(\omega)(e_n) = X'(\omega)(e_n)$ for all $n\in\N$ and all $\omega\in \Omega_1$. We may also assume that $X(\omega), X'(\omega)$ belong to $E'$ for all $\omega\in \Omega_1$. Let $e\in E$ and $e_{n_j}\to e$ in $E$. Then $X(\omega)(e) = \lim_{j\to\infty} X(\omega)(e_{n_j}) = \lim_{j\to\infty} X'(\omega)(e_{n_j}) = X'(\omega)(e)$.
Thus $X(\omega) = X'(\omega)$ for $\omega\in\Omega_1$.

The injectivity of $\bar{\kappa} = \iota_2\circ\kappa_2$ clearly implies the injectivity of $\kappa_2$.
\endproof

An example of a non-separable space for which $\bar{\kappa}$ is not injective (nor are $\kappa_0,\kappa_1,\kappa_2$) is given in Example\;\ref{rem:kappa1not1-1}.

\section{Regularity theorems}
\label{Sec:regularity}

We begin by stating the central regularization result. 

\begin{thm}\label{thm:Ito}
(a) If $E$ is a nuclear locally convex space, then $\bar{\kappa}:L^0(\Omega\!:\!E')\to \cL_{\rm c}(E\!:\!L^0(\Omega))$ is surjective.

(b) If in addition $E$ is sequentially separable, then $\bar{\kappa}$ is a bijection, i.e., the space $E$ has property ($\cR$).
\end{thm}
\emph{Proof.} The surjectivity in (a) follows from Corollary\;\ref{cor:centralresult} in the Appendix, the injectivity in (b) from Proposition\;\ref{prop:injective} above.
\endproof

The fundamental part (a) is due to It\^{o} and Nawata \cite{Ito:1984,ItoNawata:1983}. Part (b) presents a slight improvement over the uniqueness result of It\^{o} and Nawata, where almost sure uniqueness is formulated only for $\sigma$-concentrated functionals (cf. Appendix). The observation that sequential separability is an alternative hypothesis appears to be new.

\begin{rem}\label{rem:metrizablenuclear}
Every metrizable nuclear space is separable and, in particular, sequentially separable \cite[Theorem 4.4.10]{Pietsch:1972}. In general, a nuclear space need not be separable \cite[Paragraph 4.4.9]{Pietsch:1972}. Also, a separable nuclear space need not be sequentially separable; an example is $\R^I$, $I = [0,1]$, see \cite[14.4.H]{Jarchow:1981} (its nuclearity follows from \cite[Proposition 5.2.1]{Pietsch:1972}).
\end{rem}

We now turn to permanence properties, showing that surjectivity of $\bar{\kappa}$ carries over to arbitrary direct sums and to complemented subspaces, while property ($\cR$) is inherited by regular inductive limits.

\emph{Direct sums.} Let $E_\gamma$, $\gamma\in\Gamma$ be an arbitrary family of locally convex spaces and $E = \sum_{\gamma\in\Gamma}E_\gamma$ be its locally convex direct sum (i.e., topologized by the locally convex inductive limit topology with respect to the injections $\iota_\gamma:E_\gamma\to E$). Its dual is the product space $E' = \prod_{\gamma\in\Gamma} E_\gamma'$; it is endowed with the product $\sigma$-algebra with respect to the individual $\sigma$-algebras on $E_\gamma'$. We can consider the maps $\bar{\kappa}_\gamma: L^0(\Omega\!:\!E_\gamma')\to \cL_c(E_\gamma\!:\!L^0(\Omega))$ as well as
$\bar{\kappa}: L^0(\Omega\!:\!E')\to \cL_c(E\!:\!L^0(\Omega))$. Let $Y\in \cL_c(E\!:\!L^0(\Omega))$. We can consider $Y\vert E_\gamma = Y\circ\iota_\gamma$, which clearly belongs to $\cL_c(E_\gamma\!:\!L^0(\Omega))$ (since $\iota_\gamma$ is continuous).

\begin{prop}\label{prop:directsum}
Let $E = \sum_{\gamma\in\Gamma}E_\gamma$. Then $\bar{\kappa}$ is surjective if and only if all $\bar{\kappa}_\gamma$, $\gamma\in\Gamma$ are surjective,.
\end{prop}
\emph{Proof.}
The if-part:
By assumption, there are $X_\gamma \in L^0(\Omega\!:\!E_\gamma')$, $\gamma\in \Gamma$, such that $\bar{\kappa}_\gamma X_\gamma = Y|E_\gamma$. As noted earlier, we may assume that
$X_\gamma(\omega)$ is defined for all $\omega\in\Omega$. We set $X(\omega) = (X_\gamma(\omega))_{\gamma\in\Gamma} \in \prod_{\gamma\in\Gamma} E_\gamma' = E'$. Clearly, $X$ is measurable, hence belongs to $L^0(\Omega\!:\!E')$.

Recall that each $e\in E$ is of the form $e = \sum_{\rm finite} e_\gamma$ with $e_\gamma \in E_\gamma$. Then
\[
  \bar{\kappa}X(e)(\omega) = X(\omega)(e) = \sum_{\rm finite} X_\gamma(\omega)(e_\gamma) = \sum_{\rm finite} \bar{\kappa}_\gamma X_\gamma(e_\gamma)(\omega)
  =  \sum_{\rm finite} Y(e_\gamma)(\omega).
\]
Since the sum is finite, this holds on a set of full measure. Thus $\bar{\kappa}X(e) = Y(e)$ in $L^0(\Omega)$.

The only-if-part: Fix some $\gamma\in\Gamma$ and let $Y_\gamma \in \cL_c(E_\gamma\!:\!L^0(\Omega))$. Extend $Y_\gamma$ to $Y\in \cL_c(E\!:\!L^0(\Omega))$ by setting $Y \circ\iota_\delta = 0$, $\delta\neq\gamma$, and $Y\circ \iota_\gamma = Y_\gamma$. By assumption, there is $X\in L^0(\Omega\!:\!E')$ such that $\bar{\kappa} X = Y$.
Set $X_\gamma = X\circ\iota_\gamma \in L^0(\Omega\!:\!E_\gamma')$. If $f\in E_\gamma$, then $\bar{\kappa}_\gamma X_\gamma(f)(\omega) = \bar{\kappa}_\gamma X(\iota_\gamma(f))(\omega)
= X(\omega)(\iota_\gamma(f)) = Y(\iota_\gamma(f))(\omega) = Y_\gamma(f)(\omega)$. Thus $\bar{\kappa}_\gamma X_\gamma = Y_\gamma$.
\endproof

Proposition\;\ref{prop:directsum} is of interest, because uncountable direct sums of nuclear spaces are not nuclear \cite[Remark above Corollary 28.8]{MeiseVogt:1997}, yet $\bar{\kappa}$ is surjective.

\begin{rem}\label{rem:finestlctopology}
Let $\Lambda$ be any index set and $E=\R^{(\Lambda)}$. The finest locally convex topology on $E$ coincides with the direct sum topology on $\R^{(\Lambda)}$. Further, any converging sequence in $E$ is contained in a finite
dimensional subspace $F$ \cite[Spplements III(2) and V(1)]{Robertson:1980}. Using coordinates in $F$, it is quite obvious that any linear map $Y:F\to L^0(\Omega)$ is pathwise continuous. Thus
\[
    \cL(E\!:\!L^0(\Omega)) \subset \cL(E\!:\!L^0(\Omega))_{\rm path} \subset \cL(E\!:\!L^0(\Omega))_{\rm prob} \subset \cL(E\!:\!L^0(\Omega)),
\]
so all spaces are equal when $E$ is equipped with the finest locally convex topology. In particular, $\cL(E\!:\!L^0(\Omega)) = \cL_{\rm c}(E\!:\!L^0(\Omega))$. Similarly,
\[
   L^0(\Omega\!:\! E^\ast) = L^0(\Omega\!:\! E') = L^0(\Omega\!:\! E^\ast)_{\rm path}
\]
in this case. Since $\R$ is nuclear, Proposition\;\ref{prop:directsum} gives that $\kappa_0 = \kappa_1 = \bar{\kappa}$ is surjective, resulting in an alternative proof of Proposition\;\ref{prop:directsum}.
\end{rem}

\begin{cor}\label{cor:complemented}
Surjectivity of $\bar{\kappa}$ is inherited by complemented subspaces.
\end{cor}
\emph{Proof.} A complemented subspace $E_1$ of a locally convex space possesses, by definition, a topological complement $E_2$, so that $E$ equals the direct sum $E_1\oplus E_2$. Apply Proposition\;\ref{prop:directsum}. \endproof

\emph{Countable inductive limits.} Let $E_1\subset E_2 \subset \ldots \subset E_n \subset \ldots $ be a sequence of locally convex spaces with continuous inclusion maps. Let the linear space $E$ be their union. If $E$ carries the finest locally convex topology making all inclusions $E_n\subset E$ continuous, then it is called the \emph{inductive limit} of the spaces $E_n$. The inductive limit is called \emph{regular} if each bounded subset of $E$ is contained in some $E_n$ and bounded there
\cite[Section 4.5]{Jarchow:1981}. As above, we can consider the mappings
$\bar{\kappa}_n: L^0(\Omega\!:\!E_n')\to \cL_c(E_n\!:\!L^0(\Omega))$ as well as $\bar{\kappa}: L^0(\Omega\!:\!E')\to \cL_c(E\!:\!L^0(\Omega))$.

\begin{prop}\label{prop:indlim}
Let $E$ be the regular inductive limit of the spaces $E_n$. If all $\bar{\kappa}_n$, $n\in\N$ are bijective, then $\bar{\kappa}$ is bijective.
\end{prop}
\emph{Proof.}
Let $Y\in\cL_c(E\!:\!L^0(\Omega))$. This means that if $e_k\to e$ in $E$, then $Y(e_k)\to Y(e)$ in $L^0(\Omega)$. Fix $n\in\N$. If $e_k\to e$ in $E_n$ then $e_k\to e$ in $E$.
Thus $Y|E_n \in \cL_c(E_n\!:\!L^0(\Omega))$. By assumption, there exists $X_n\in L^0(\Omega\!:\!E_n')$ such that $\bar{\kappa}_n X_n = Y|E_n$. Further, if $X_n'\in L^0(\Omega\!:\!E_n')$ also satisfies $\bar{\kappa}_n X_n' = Y|E_n$, then $X_n = X_n'$ almost surely.

Next, observe that for all $e\in E_n$, $\bar{\kappa}_{n+1}X_{n+1}(e) = Y(e)$ in $L^0(\Omega)$. This implies that
$\bar{\kappa}_{n}(X_{n+1}|E_n)(e) = Y(e) = \bar{\kappa}_{n}X_{n}(e)$. Thus $\bar{\kappa}_{n}(X_{n+1}|E_n) = \bar{\kappa}_{n}X_{n} = Y|E_n$.
The injectivity of $\bar{\kappa}_{n}$ implies that there is a set $\Omega_n\subset\Omega$ of full measure such that
$X_{n+1}(\omega)|E_n = X_n(\omega)$ for all $\omega\in \Omega_n$.

Let $\Omega' = \bigcap_{n=1}^\infty\Omega_n$. Then $P(\Omega') = 1$ and we have obtained a sequence of functionals $X_n\in L^0(\Omega\!:\!E_n')$
such that $X_{n+m}(\omega)(e) = X_n(\omega)(e)$ for all $m,n\in\N$, $e\in E_n$ and $\omega\in\Omega'$. For $\omega\in\Omega'$, the map
\[
   X(\omega):E = \bigcup_{n=1}^\infty E_n\to \R
\]
given by $X(\omega)(e) = X_{n}(\omega)(e)$ for $e\in E_n$ is well-defined. It is clearly linear. Let $e_k\to e$ in $E$. By regularity, there is $m\in\N$ such that all $e_k$ and $e$ belong to $E_m$. The topology of $E_m$ is finer than the one induced by $E$ on $E_m$, hence $e_k\to e$ in $E_m$, and so $X(\omega)(e_k) = X_m(\omega)(e_k) \to X_m(\omega)(e) = X(\omega)(e)$. Thus $X(\omega):E\to \R$ is sequentially continuous and hence belongs to $E'$, the sequential dual of $E$.
Since $\omega\to X(\omega)(e)$ is measurable for all $e\in E$, $X$ is measurable when $E'$ is endowed with the Kolmogorov $\sigma$-algebra. Summarizing, we have obtained
$X\in L^0(\Omega\!:\!E')$ such that $\bar{\kappa}X = Y$, and $\bar{\kappa}$ is seen to be surjective.

Assume that $X'\in L^0(\Omega\!:\!E')$ such that $\bar{\kappa}X' = Y$. This means in particular for $e\in E_n$ that $X'(\omega)(e) = Y(e)(\omega) = X_n(\omega)(e)$ on a set of full measure; thus $\bar{\kappa}_n(X'|E_n) = \bar{\kappa}_n X_n$. By the assumed injectivity of $\bar{\kappa}_n$, there is a set $\Omega_n'\subset\Omega$ of full measure such that
$X'(\omega)|E_n = X_n(\omega)$ for all $\omega\in \Omega_n'$. Let $\Omega'' = \bigcap_{n=1}^\infty\Omega_n'$. Then $P(\Omega'') = 1$ and $X'(\omega) = X(\omega)$ for $\omega\in\Omega''$, so $\kappa$ is seen to be injective.
\endproof

\begin{ex}\label{ex:distributionspaces}
The following spaces of smooth functions have property ($\cR$):
\[
  \cS(\R^n),\ \cE(\R^n),\ \cD_K(\R^n),\ \cD(\R^n)
\]
where $K\subset \R^n$ is compact. The first three spaces are metrizable nuclear spaces (see e.g. \cite[Section 51]{Treves:1967}, \cite[Section 6.2]{Pietsch:1972}). 
The space $\cD(\R^n)$ is the strict inductive limit of the spaces $\cD_K(\R^n)$, so one can either use Proposition\;\ref{prop:indlim} 
or invoke its nuclearity \cite[Section 1.5]{Ito:1984} and sequential separability.
\end{ex}

\section{Counterexamples and mapping properties}
\label{Sec:counter}

{\it Counterexamples to injectivity.} The meaning of the maps $\kappa, \kappa_0, \kappa_1, \kappa_2, \bar{\kappa}$ and $\iota, \iota_0, \iota_1, \iota_2$ can be recalled from the diagram in Section\;\ref{Sec:contlinrandfun}.
\begin{ex}\label{rem:kappanot1-1}
The map $\kappa$ is not necessarily injective. Indeed, take $\Omega = \R$ with any probability measure which has a continuous, everywhere positive density, and take $E = \R$.
Let
\[
   X(\omega)(e) = \left\{\begin{array}{ll}
                      1, & e = \omega,\\
                      0, & e \neq \omega.
                      \end{array}\right.
\]
Then $\{\omega\in\Omega:\forall e \in E, X(\omega)(e)=0\} = \emptyset$, so $[X] \neq [0]$. But $[\kappa X(e)] = 0$ in $L^0(\Omega)$ for all $e\in E$,
because $\{\omega\in\Omega: X(\omega)(e) \neq 0\} = \{e\}$, and $P(\{e\}) =0$.
\end{ex}
The process $X$ in Example\;\ref{rem:kappanot1-1} is a version of the zero operator $Y \equiv 0$ which belongs, in particular, to $\cL_{\rm c}(E\!:\!L^0(\Omega))$. Note that $X$ is not only nonzero, but also
nonlinear, since
\[
   \{\omega \in\Omega: \forall e,f, X(\omega)(e + f) = X(\omega)(e) +  X(\omega)(f)\} = \emptyset.
\]
Thus $X\in L^0(\Omega\!:\!\R^E)$ is a version of $Y\equiv 0$ which does not belong to $L^0(\Omega\!:\!E^\ast)$.
\begin{ex}\label{rem:kappa0not1-1}
A generalization of Example\;\ref{rem:kappanot1-1} shows that again $\kappa_0$ need not be injective, in general. To see this, take $\Omega = \R$ with a probability measure as in
Example\;\ref{rem:kappanot1-1} and set $E = \R^{(\R)}$. Equip $E^\ast =\R^{\R}$ with the Kolomogorov $\sigma$-algebra $\cB_K(\R^{\R})$. Define
\[
   X \in L^0(\Omega\!:\!E^\ast)\ {\rm by\ } X(\omega)(e) = \pi_\omega(e) = e(\omega),\ e\in E.
\]
Then $\{\omega\in\Omega:\forall e \in E, X(\omega)(e)=0\} = \emptyset$, so $[X] \neq [0]$. But $[\kappa_0 X(e)] = 0$ in $L^0(\Omega)$
because $\{\omega\in\Omega: X(\omega)(e) \neq 0\}$ is finite for all $e\in E$, hence has zero probability.
\end{ex}

\begin{ex}\label{rem:kappa1not1-1}
If $E$ in Example\;\ref{rem:kappa0not1-1} is given the finest locally convex topology (see Remark\;\ref{rem:finestlctopology}), then $E' = E^\ast$. The process $X$ in Example\;\ref{rem:kappa0not1-1}
belongs to $L^0(\Omega\!:\!E^\ast) = L^0(\Omega\!:\!E') = L^0(\Omega\!:\!E^\ast)_{\rm prob}$, and $X\neq 0$ as a member of those spaces. On the other hand,
$\kappa_0 X(e) = \kappa_1 X(e) = \kappa_2 X(e) = \bar{\kappa} X(e) = 0$ in $L^0(\Omega)$ for all $e\in E$. This shows that $\kappa_1$, $\kappa_2$ and
$\bar{\kappa}$ are not injective in this situation.

Note that $\R^{(\R)}$ is not separable \cite[14.4.H]{Jarchow:1981}, so this example shows that the hypothesis of Proposition\;\ref{prop:injective} cannot be dropped, in general.
\end{ex}

{\it Counterexamples to surjectivity.} For completeness, we state possible lack of surjectivity of $\iota$ and $\iota_0$. It trivially rests on the fact that there are locally convex spaces $E$ such that
$E^\ast \neq E'$.
\begin{ex}\label{rem:iotanotsurjective}
There is a probability space $(\Omega,\Sigma,P)$ and a locally convex space $E$ such that $\iota$ and $\iota_0$ are not surjective.

In fact, it suffices to take any space $E$ such that $E'\neq E^\ast$ and a one-point probability space $\Omega$. If $X\in E^\ast\setminus E'$, there is $e_n\to e$ such that
$X(e_n)\not\to X(e)$. There is $\eps > 0$ such that $|X(e_n)- X(e)| > \eps$ for infinitely many $n$. As every event on a one-point probability space, this event has probability 1. Thus
$x\in \cL^0(\Omega\!:\!E^\ast)\setminus \cL^0(\Omega\!:\!E^\ast)_{\rm prob}$, so $\iota$ is not surjective. The same reasoning applies to $\iota_0$, considering $X$ as a map $E\to L^0(\Omega) = \R$.
\end{ex}

The next examples show that Theorem\;\ref{thm:Ito}(a) does not hold in general: without the nuclearity assumption, $\bar{\kappa}$ need not be surjective. It may even happen that $L^2$-continuous elements of $L^0(\Omega\!:\!E^\ast)$
do not have a version in $E'$.

\begin{ex}\label{ex:noversioninEprime}
We construct a Gaussian probability space $(\Omega,\Sigma,P)$, a normed space $E$ and a linear random functional $X\in L^0(\Omega\!:\!E^\ast)_{\rm prob}$, $\bar{\kappa} X\in \cL_{\rm c}(E\!:\!L^2(\Omega))$ such that
$P(\{\omega\in\Omega: X(\omega)\not\in E'\}) = 1$.
\end{ex}

Indeed, let $\Omega_0 = (\R,\cB(\R),P_0)$ where $P_0$ is the law of a standard normal distribution $\cN(0,1)$. We let $\Omega$ be the product space $\Omega_0^\N$, more precisely,
$\Omega = (\R^\N, \cB_K(R^\N), \otimes P_0)$ where $P =\otimes P_0$ is the product probability.

Next, let $E = \R^{(\N)}$ eqipped with the $\ell^2$-norm (actually any $\ell^p$ with $1\leq p \leq 2$ would do). Since $E$ is dense in $\ell^2$, we have $E' = \ell^2$, while $E^\ast = \R^\N$.
The elements of $E$ are of the form
\[
   e = \sum_{\rm finite}\alpha_j e_j
\]
with $\alpha_j\in\R$ and $e_j$ the $j$th unit vector.

Define the linear random functional $X$ by $X(\omega)(e) = \langle\omega,e\rangle$. Clearly, $X(\omega)\in E^\ast$ for all $\omega$ and
\[
   \langle\omega,e\rangle = \sum_{\rm finite}\omega_j\alpha_j \quad {\rm when}\quad e = \sum_{\rm finite}\alpha_j e_j.
\]
As a finite sum of measurable functions, $\omega\to X(\omega)(e)$ is measurable for all $e\in E$, so $X\in L^0(\Omega\!:\!E^\ast)$.

\emph{The probability distribution of $X(\cdot)(e)$.} Clearly, it is Gaussian with mean zero, $\rE\big(X(\cdot)(e)\big) = 0$ and variance
\[
\rV\big(X(\cdot)(e)\big) = \rE\big(X(\cdot)(e)^2\big) = \rE\Big(\sum_{\rm finite}\omega_j\alpha_j\sum_{\rm finite}\omega_k\alpha_k\Big) = \sum_{\rm finite}\alpha_j^2 = \|e\|_{\ell^2}^2
\]
since $\rE(\omega_j\omega_k) = \delta_{jk}$.

Let $e^{(n)}, e \in E$ and $e^{(n)} \to e$ in $\ell^2$. Then
\[
   \rE\Big(\big(X(\cdot)(e^{(n)}) - X(\cdot)(e)\big)^2\Big) = \rE\Big(\big(X(\cdot)(e^{(n)} - e)\big)^2\Big) = \|e^{(n)} - e\|_{\ell^2} \to 0,
\]
thus $\kappa X \in \cL_{\rm c}(E\!:\!L^2(\Omega))$.

Finally, we show that for almost all $\omega\in\Omega$, $X(\omega)\not\in E'$. Note that
\[
   \{\omega\in\Omega:X(\omega)\not\in E'\} = \{\omega\in\Omega: \exists e^{(n)}\to e, X(\omega)(e^{(n)}) \not\to X(\omega)(e)\}
\]
Let $e^{(n)} = \alpha_n e_n$ with $\alpha_n\to 0$ to be specified. Then $e^{(n)} \to 0 \equiv e$ in $\ell^2$. Let
\[
   A_n = \{\omega\in\Omega: X(\omega)(e^{(n)}) \geq 1\} = \{\omega\in\Omega: \omega_n\alpha_n\geq 1\} = \{\omega\in\Omega: \omega_n\geq\tfrac{1}{\alpha_n}\}.
\]
Let
\[
   A = \{\omega\in\Omega: \omega\in A_n\ {\rm for\ infinitely\ many\ }n\}.
\]
Note that the $A_n$ are independent. With this choice of $e^{(n)}$ and $e\equiv 0$,
\[
   A \subset \{\omega\in\Omega, X(\omega)(e^{(n)}) \not\to 0\} \subset \{\omega\in\Omega, X(\omega)\not\in E'\}.
\]
We will show that $\sum_{n=1}^\infty P(A_n) = \infty$. The Borel-Cantelli-Lemma then implies that $P(A) = 1$. But
\[
   P(A_n) = \frac{1}{\sqrt{2\pi}}\int_{1/\alpha_n}^\infty \ee^{-y^2/2}\dd y \geq \frac{1}{\sqrt{2\pi}}\int_{1/\alpha_n}^\infty y\ee^{-y^2}\dd y
       = \frac{1}{2\sqrt{2\pi}}\ee^{-1/\alpha_n^2}.
\]
Choosing $\alpha_n = (\log n)^{-1/2}$, we have $P(A_n) \geq n/2\sqrt{2\pi}$. With this choice,
$P(A) = 1 \leq P\big(\{\omega\in\Omega, X(\omega)\not\in E'\}\big)$.

\begin{ex}\label{rem:iotasnotonto}
Example\;\ref{ex:noversioninEprime} can be further exploited to show that the maps $\iota_1$ and $\iota_2$ need not be surjective.
It follows from its proof that the process $X$ is $L^2$-continuous: if $e^{(n)}\to e$ in $E$, then $X(\cdot)(e^{(n)}-e)\to 0$ in $L^2(\Omega)$.
This implies that $X\in \cL(\Omega\!:\!E^\ast)_{\rm prob}$, but it does not belong to $\cL(\Omega\!:\!E') = \cL(\Omega\!:\!E^\ast)_{\rm path}$.
Thus $\iota_1$ is not surjective.

The statement $P\big(\{\omega\in\Omega, X(\omega)\not\in E'\}\big) = 1$ can also be read as
\[
   P\big(\{\omega\in\Omega: \exists e^{(n)}\to e, \kappa_2X(e^{(n)})(\omega) \not\to \kappa_2X(e)(\omega)\}\big) = 1,
\]
showing that $\kappa_2 X$ does not belong to $\cL(E\!:\!L^0(\Omega))_{\rm path}$, while clearly $\kappa_1 X$ belongs to $\cL(E\!:\!L^0(\Omega))_{\rm prob}$.
Thus $\iota_2$ is not surjective either.
\end{ex}

\begin{ex}\label{ex:Lp}
Let $\Omega = (0,1)$ with the Borel $\sigma$-algebra $\Sigma = \cB(\R)$ and the uniform probability measure $P$. Let $E = L^p(0,1)$, where $1\leq p < \infty$. Define $Y\in \cL_c(E\!:\!L^0(\Omega))$ by $Y(e) = e$. Then there does not exist $X\in L^0(\Omega\!:\!E')$ such that
$\bar{\kappa}X = Y$. Indeed, such an $X$ would have to satisfy $X(\omega)(e) = e(\omega)$ on some set of full measure $\Omega_e$, possibly depending on $e$. There are two reasons why this cannot happen. First, let $(e_n)_{n\in\N}$ be any convergent sequence in $L^p(0,1)$, say with limit $e$.
We have $X(\omega)(e_n) = e_n(\omega)$ on the set $\Omega' = \bigcap_{n=1}^\infty \Omega_{e_n} \cap \Omega_e$, which has full measure. If $X(\omega)$ were a member of $E'$, then $X(\omega)(e_n) \to X(\omega)(e)$ for $\omega\in\Omega'$.
But this means that $e_n(\omega) \to e(\omega)$ almost everywhere. The argument implies that every $L^p$-converging sequence converges almost everywhere, which is not true.

Second, one could also argue that $X(\omega)$ is a point-evaluation functional, which cannot belong to $E' = L^q(0,1)$, $\frac1p + \frac1q = 1$. This is excluded because $X(\omega)$ restricted to $C[0,1]$ is the evaluation functional (Dirac measure) at point $\omega$. It is wellknown that the evaluation functional on $C[0,1]$ does not have an $L^p$-continuous extension to $L^p(0,1)$.
\end{ex}
Example\;\ref{ex:Lp} is another instance where $\bar{\kappa}$ is not surjective. We now produce an example in which even $\kappa_2$ is not surjective, thus even pathweise continuous linear random functionals need not have a version in the dual.

\begin{ex}\label{ex:Linfty}
Take again $\Omega = (0,1)$, $E = L^\infty(0,1)$ and $Y\in \cL_c(E\!:\!L^0(\Omega))$ given by $Y(e) = e$. Clearly, $Y$ belongs to $\cL(E\!:\!L^0(\Omega))_{\rm path}$. There does not exist $X\in L^0(\Omega\!:\!E')$ such that $\bar{\kappa}X = Y$. Indeed, we would again have that $X(\omega)(e) = e(\omega)$, so $X(\omega)$, restricted to $C[0,1]$, would simply be the Dirac measure at $\omega$. According to \cite[Theorem IV.8.16]{DunfordSchwartz:1958}, the dual of $L^\infty(0,1)$ is equal to the space of bounded additive set functions on the Lebesgue extension of $\Sigma$, which vanish on sets of Lebesgue measure zero. The latter is obviously false for the Dirac measure, which hence cannot belong to the dual of $L^\infty(0,1)$.
\end{ex}

We can summarize the mapping properties as follow.

(a) As was shown in Sections\;\ref{Sec:randfun} and \ref{Sec:linrandfun}, the maps $\kappa$ and $\kappa_0$ are always surjective.

(b) The maps $\kappa$, $\kappa_0$, $\kappa_1, \kappa_2$ and $\bar{\kappa}$ need not be injective (Examples\;\ref{rem:kappanot1-1}, \ref{rem:kappa0not1-1}, \ref{rem:kappa1not1-1}).

(c) The maps $\iota$ and $\iota_0$ need not be surjective (Example\;\ref{rem:iotanotsurjective}).

(d) The maps $\bar{\kappa}$, $\iota_1$, $\iota_2$ and $\kappa_2$ are not necessarily surjective (Examples\;\ref{ex:noversioninEprime} and \ref{ex:Lp}, Example\;\ref{rem:iotasnotonto}, Example\;\ref{ex:Linfty}).

It remains open whether $\kappa_1$ is surjective, in general. If a locally convex space has property ($\cR$), then $\bar{\kappa}$ as well as $\kappa_2$ and $\iota_2$ are bijective, and $\kappa_1$ is surjective
(Remark\;\ref{rem:consequencesOfR}).

\appendix
\section{Appendix: Regularity and nuclear spaces}
\label{Sec:appendix}

The present paper relies on the theorem that for every nuclear locally convex space the map $\bar{\kappa}$ is surjective. In a slightly more general form, the theorem has been proved by It\^{o} and Nawata \cite{ItoNawata:1983} and in It\^{o} \cite{Ito:1984}. A somewhat less general version, simplifying It\^{o}'s proof, is due to Walsh \cite{Walsh:1984}.
The purpose of this appendix is to record It\^{o}'s version \cite{Ito:1984} of the theorem for reference in the paper, and to compare it with the version of Walsh \cite{Walsh:1984}.
The proofs of this result are easily readable both in \cite{Ito:1984} and \cite{Walsh:1984} and will not be reproduced here.

Following It\^{o}'s presentation \cite[Section 1.2]{Ito:1984}, a topology $\tau$ on a linear space $E$ is called \emph{multi-Hilbertian}, if it is generated by a directed family $\cP$ of separable Hilbertian seminorms. If $\cP$ is countable, $(E,\tau)$ is referred to as \emph{countably Hilbertian}.

Let $p, q$ be two separable Hilbertian seminorms on $E$. Then $p$ is called \emph{Hilbert-Schmidt bounded} by $q$ if $\sum_{n=1}^\infty p(e_n)^2$ is finite for some (and hence any \cite[Remark 1.1.2]{Ito:1984}) orthonormal basis $(e_n)_{n\in \N}$ in $(E,q)$.
Let $\tau_1, \tau_2$ be multi-Hilbertian topologies on $E$. Then $\tau_1$ is termed \emph{Hilbert-Schmidt weaker} than $\tau_2$, if every $\tau_1$-continuous (separable Hilbertian) seminorm $p$ is Hilbert-Schmidt bounded by some $\tau_2$-continuous (separable Hilbertian) seminorm $q$. Given a multi-Hilbertian topology $\tau$ on $E$,
the strongest multi-Hilbertian topology Hilbert-Schmidt weaker than $\tau$ is called the \emph{Kolmogorov $I$-topology of $\tau$} and denoted by $I(\tau)$. It is generated by the family of separable Hilbertian seminorms which are Hilbert-Schmidt weaker than some $\tau$-continuous (separable Hilbertian) seminorm.

For the notion of a \emph{nuclear locally convex space} we follow the classical sources \cite[Chap. II \S 2, n$^{\rm o}$ 1]{Grothendieck:1955}, \cite[Sect. 4.1]{Pietsch:1972}, see also \cite{Jarchow:1981,Schaefer:1980,Treves:1967}.

\begin{prop}\label{prop:equivalenceOfNuclearity}
A locally convex space $(E,\tau)$ is nuclear if and only if it is multi-Hilbertian and $I(\tau) = \tau$.
\end{prop}
\emph{Proof.}
Let $\cP$ be a directed fundamental system of seminorms defining the topology $\tau$. For $p\in\cP$, let $E_p$ be the completion of the normed space $E/p^{-1}(0)$. Recall that the topology $\tau$ on $E$ is equal to the projective topology with respect to the canonical maps $\iota^p:E\to E_p$, $p\in\cP$ \cite[Remark 24.5(a)]{MeiseVogt:1997}. For $p \leq q$, the identity map induces a continuous linear map $\iota_q^p: E_q\to E_p$ such that $\iota_q^p \circ \iota^q = \iota^p$.

According to \cite[Proposition 28.4]{MeiseVogt:1997}, a locally convex space is nuclear if and only if it has a fundamental system $\cP$ of Hilbertian seminorms and for every continuous Hilbertian seminorm $p$ there exists a continuous Hilbertian seminorm $q$, $p\geq q$ such that
$\iota_q^p: E_q\to E_p$ is a Hilbert-Schmidt operator. What is more, using \cite[Lemma 29.7]{MeiseVogt:1997} and its proof together with \cite[Remark 24.5(b)]{MeiseVogt:1997}, for each continuous Hilbertian seminorm $p$, there is a countable orthonormal basis in $E_p'$. In particular, $E_p'$ and consequently $E_p$ are separable. 
Thus every nuclear space is multi-Hilbertian. 

It is easy to see that $p$ is Hilbert-Schmidt bounded by $q$ if and only if the map $\iota_q^p: E_q\to E_p$ is a Hilbert-Schmidt operator. (The essential argument can be found e.g. in \cite[Section 2.5]{Pietsch:1972}. To go over to the quotient spaces, one may use \cite[Remark 12.5(a)]{MeiseVogt:1997} together with \cite[Theorem 1.1.1(i)]{Ito:1984}.) 

It follows that for a multi-Hilbertian space, the assertion that $I(\tau)=\tau$ is equivalent with nuclearity.
\endproof

In order to formulate the central \emph{regularization theorem}, we recall some terminology from \cite[Section 2.3]{Ito:1984}. As in the present paper, the elements of $\cL(E\!:\!L^0(\Omega))$ are called \emph{linear random functionals} on $E$.
Let $E$ be equipped with a multi-Hilbertian topology $\tau$. The elements of $\cL_{\rm c}(E\!:\!L^0(\Omega))$ are referred to as \emph{$\tau$-continuous linear random functionals} in 
\cite[(2) in Section 2.3]{Ito:1984}.
An $E_\tau'$-valued random variable $X$ (i.e., an element of $L^0(\Omega\!:\!E_\tau')$) is called \emph{$\sigma$-concentrated} if there exists a countably Hilbertian topology $\theta$ coarser than $\tau$ such that $P(X\in E_\theta') = 1$, \cite[Definition 2.3.1]{Ito:1984}. Let $Y\in \cL_{\rm c}(E\!:\!L^0(\Omega))$. Then an $E_\tau'$-valued random $X$ is called a \emph{$\tau$-regular version} of $Y$ if (i) $X$ is $\sigma$-concentrated and (ii) for every $e\in E$, $X(e) = Y(e)$ almost surely (\cite[Definition 2.3.1]{Ito:1984}). In the language of Section\;\ref{Sec:contlinrandfun}, (ii) means that $\bar{\kappa}X = Y$. Collecting Theorems 2.3.1, 2.3.2 and 2.3.3 from \cite{Ito:1984} gives the central result:
\begin{thm}(The regularization theorem \cite{Ito:1984,ItoNawata:1983})
\label{thm:centralresult}
A linear random functional $Y$ on $E$ has a $\tau$-regular version if and only if it is $I(\tau)$-continuous. It is unique in the sense that if $X$ and $X'$ are $\tau$-regular versions of $Y$, then $X=X'$ almost surely.
\end{thm}
\begin{cor}\label{cor:centralresult}
If $E$ is a nuclear space, then $\bar{\kappa}:L^0(\Omega\!:\!E')\to \cL_{\rm c}(E\!:\!L^0(\Omega))$ is surjective.
\end{cor}
\emph{Proof.} It follows by combining Theorem\;\ref{thm:centralresult} with Proposition\;\ref{prop:equivalenceOfNuclearity}.
\endproof

Theorem\;\ref{thm:centralresult} states almost sure uniqueness only for $\tau$-regular, hence $\sigma$-con\-cen\-trated versions. It is not evident whether this gives injectivity of $\bar{\kappa}$ in general.

There is an important special case of countably Hilbertian nuclear spaces. Suppose that the topology of a countably Hilbertian nuclear space $E$ is generated by a directed
sequence of separable Hilbertian norms $p_1 \leq p_2 \leq\ldots\leq p_n \leq \ldots$. In addition, suppose that the norms are compatible \cite[Chapter I.2.2]{GelfandShilov:1968}. 

The Hilbert spaces $E_m = E_{p_m}$ are defined as above (as the completion of $E$ with respect to $p_m$).
Due to compatibility, the maps $\iota_n^m:E_n\to E_m$ ($n \geq m$) are injective. One obtains the nested sequence of Hilbert spaces
$E_1 \supset E_2 \supset \ldots  E_n \supset \ldots $ with continuous inclusion maps. The space $E$ is equipped with the corresponding projective topology.
The space $E$ is complete if and only if it coincides with the intersection
$
   E = \bigcap_{n=1}^\infty E_n
$,
see \cite[Theorem in Chapter I.3.2]{GelfandShilov:1968} (and, in particular, a Fr\'{e}chet space). 
By the argument given in  the proof of Propositon\;\ref{prop:equivalenceOfNuclearity}, each $E_n$ is separable; for $E$, separability holds because every metrizable nuclear space is separable \cite[Theorem 4.4.10]{Pietsch:1972}. (For the separability results, one may also refer to \cite[Chapter I.6.5]{GelfandShilov:1968} and
\cite[Chapter I.3.4]{GelfandVilenkin:1964}.) In particular, a nuclear space $E= \bigcap_{n=1}^\infty E_n$ as above is countably Hilbertian and a separable Fr\'{e}chet space. This class of spaces has been considered in Gel'fand-Vilenkin \cite[Chapter I.3.2]{GelfandVilenkin:1964}. 

The proof of Walsh \cite[Theorem 4.1 and Corollary 4.2]{Walsh:1984} shows that $\bar{\kappa}:L^0(\Omega\!:\!E')\to \cL_{\rm c}(E\!:\!L^0(\Omega))$ is surjective for this narrower class.

\end{document}